\documentclass[reqno]{amsart}
\usepackage{amssymb,eucal,graphics}

\newcommand{\jz}{\ensuremath{\langle\vee,0\rangle}}

\newcommand{\two}{\mathbf{2}}

\newcommand{\bp}{\mathbin{\square}}
\newcommand{\ltp}{\mathbin{\boxtimes}}

\newcommand{\ootimes}{\mathbin{\ol{\otimes}}}

\def\dnw{\mathop{\downarrow}}
\def\upw{\mathop{\uparrow}}

\newcommand{\DR}{\mathbin{D}}
\newcommand{\dr}{\mathbin{\delta}}
\newcommand{\DD}{\Delta}
\newcommand{\lDD}{\overline{\DD}}

\newcommand{\tr}{\vartriangleleft}
\newcommand{\ntr}{\ntriangleleft}
\newcommand{\utr}{\trianglelefteq}
\newcommand{\nutr}{\ntrianglelefteq}
\newcommand{\dtr}{\mathbin{\vartriangleleft\kern-10pt
{\lower3pt\hbox{$\scriptscriptstyle\ne$}}\kern3pt}}
\newcommand{\trc}{\tr^{\mathrm{c}}}

\newcommand{\jzs}{\jz-se\-mi\-lattice}
\newcommand{\cm}{com\-mu\-ta\-tive mo\-no\-id}
\newcommand{\jirr}{join-ir\-re\-duc\-i\-ble}

\newcommand{\pind}{pseu\-do-in\-de\-com\-pos\-a\-ble}
\newcommand{\jsd}{join-sem\-i\-dis\-trib\-u\-tive}
\newcommand{\jsdy}{join-sem\-i\-dis\-trib\-u\-tiv\-i\-ty}
\newcommand{\msd}{meet-sem\-i\-dis\-trib\-u\-tive}

\newcommand{\set}[1]{\{#1\}}
\newcommand{\setm}[2]{\{#1\mid#2\}}

\newcommand{\lx}{\boldsymbol{x}}
\newcommand{\ly}{\boldsymbol{y}}

\newcommand{\la}{\boldsymbol{a}}
\newcommand{\lb}{\boldsymbol{b}}

\newcommand{\olp}{\overline{p}}
\newcommand{\olq}{\overline{q}}

\newcommand{\otp}{\tilde{p}}
\newcommand{\otq}{\tilde{q}}

\newcommand{\es}{\varnothing}
\newcommand{\into}{\hookrightarrow}

\newcommand{\eps}{\varepsilon}

\newcommand{\CM}{\mathcal{M}}

\newcommand{\ol}[1]{\overline{#1}}
\newcommand{\FL}{\mathrm{F}_{\mathbf{L}}}

\DeclareMathOperator{\Max}{Max}

\newcommand{\diag}[1]{{#1}^{[2]}}

\newcommand{\Co}{\mathbf{Co}}
\newcommand{\CB}{\mathbf{CB}}

\DeclareMathOperator{\J}{J}
\DeclareMathOperator{\JC}{J^c}

\DeclareMathOperator{\Dim}{Dim}
\DeclareMathOperator{\Conc}{Con_c}
\DeclareMathOperator{\Con}{Con}
\DeclareMathOperator{\Ant}{Ant}

\newcommand{\EE}{\mathbf{E}}
\newcommand{\FF}{\mathbf{F}}

\DeclareMathOperator{\R}{R}

\newcommand{\ZZ}{\mathbb{Z}}
\newcommand{\RR}{\mathbb{R}}
\newcommand{\ZZb}{\overline{\ZZ}^+}

\newcommand{\case}[1]{\smallskip\noindent{\textbf{\textit{Case}~{#1}.}}} 

\numberwithin{equation}{section}

\theoremstyle{plain}

\newtheorem{lemma}{Lemma}[section]
\newtheorem{theorem}[lemma]{Theorem}
\newtheorem{proposition}[lemma]{Proposition}
\newtheorem{corollary}[lemma]{Corollary}
\newtheorem*{sclaim}{Claim}
\newtheorem{claim}{Claim}

\newtheorem*{stat}{\name}
\newcommand{\name}{testing}

\theoremstyle{definition}
\newtheorem{definition}[lemma]{Definition}
\newtheorem{problem}{Problem}
\newtheorem{example}[lemma]{Example}
\newtheorem*{notation}{Notation}

\theoremstyle{remark}
\newtheorem{remark}[lemma]{Remark}
\newtheorem*{note}{Note}

\newenvironment{all}[1]{\renewcommand{\name}{#1}\begin{stat}} 
{\end{stat}}

\newcommand{\qedc}{{\qed}~{\rm Claim~{\theclaim}.}}
\newcommand{\qedsc}{{\qed}~{\rm Claim.}} 

\newenvironment{cproof} {\begin{proof}[Proof of Claim.]}
{\qedc\renewcommand{\qed}{}\end{proof}}

\newenvironment{scproof} {\begin{proof}[Proof of Claim.]}
{\qedsc\renewcommand{\qed}{}\end{proof}}

\begin{document}

\title[Dimension theory of lattices]
{{}From join-irreducibles to dimension theory\\
for lattices with chain conditions}

\author[F.~Wehrung]{Friedrich Wehrung}
\address{CNRS, UMR 6139\\
D\'epartement de Math\'ematiques\\
Universit\'e de Caen\\ 14032 Caen Cedex\\
France}
\email{wehrung@math.unicaen.fr}
\urladdr{http://www.math.unicaen.fr/\~{}wehrung} 


\date{\today}

\begin{abstract}
For a finite lattice $L$, the congruence lattice
$\Con L$ of $L$ can be easily computed from the partially ordered
set $\J(L)$ of \jirr\ elements of $L$ and the
\emph{join-dependency relation} $\DR_L$ on $\J(L)$. We establish a
similar version of this result for the \emph{dimension monoid} $\Dim L$
of $L$, a natural precursor of $\Con L$. For $L$ \jsd, this result
takes the following form:

\begin{all}{Theorem 1}
Let $L$ be a finite \jsd\ lattice. Then
$\Dim L$ is isomorphic to the \cm\ defined by generators $\DD(p)$,
for $p\in\nobreak\J(L)$, and relations
 \[
 \DD(p)+\DD(q)=\DD(q),\text{ for all }p,\,q\in\J(L)\text{ such
 that } p\DR_Lq.
 \]
\end{all}

As a consequence of this, we obtain the following results: 

\begin{all}{Theorem 2}
Let $L$ be a finite \jsd\ lattice. Then $L$
is a lower bounded homomorphic image of a free lattice if{f} $\Dim L$
is \emph{strongly separative}, if{f} it satisfies the axiom 
 \[
 (\forall x)(2x=x\Rightarrow x=0).
 \]
\end{all}

\begin{all}{Theorem 3}
Let $A$ and $B$ be finite \jsd\ lattices.
Then the \emph{box product} $A\bp B$ of $A$ and $B$ is \jsd, and
the following isomorphism holds: 
 \[
 \Dim(A\bp B)\cong\Dim A\otimes\Dim B.
 \]
\end{all}

\end{abstract}

\maketitle

\section{Introduction}\label{S:Intro}

The classical dimension theory of complemented modular lattices, and,
more particularly, the \emph{continuous geometries} (\emph{i.e.},
complete, upper continuous, and lower continuous complemented modular
lattices), originates in work by von~Neumann, see J. von
Neumann \cite{Neum60} or F. Maeda \cite{FMae58}. It has been
established that the von~Neumann dimension in a
continuous geometry is a particular case of a notion of dimension
defined for \emph{any} lattice. This
dimension is materialized by the so-called \emph{dimension monoid}
$\Dim L$ of a lattice $L$, see F. Wehrung \cite{WDim}.

The dimension monoid of $L$ is generated by ``distances'' $\DD(x,y)$,
for $x\leq y$ in~$L$. The compact congruence semilattice $\Conc L$ of
$L$ is the maximal semilattice quotient of $\Dim L$, and the
generator $\DD(x,y)$ is sent, \emph{via} the canonical projection, to
the principal congruence $\Theta(x,y)$.
For an irreducible continuous geometry
$L$, $\Dim L$ is isomorphic either to the chain $\ZZ^+$ of natural
numbers or to the chain $\RR^+$ of nonnegative real numbers. For a
\emph{reducible} continuous geometry $L$, the dimension monoid is the
positive cone of a Dedekind complete lattice-ordered group (see
T.~Iwamura~\cite{Iwam44}), hence, if $L$ is a bounded lattice,
dimensionality in $L$ is described by a \emph{family} of real-valued
dimension functions.

If $L$ is \emph{modular}, then dimensionality in $L$ is related to
\emph{perspectivity}, for example,
 \[
 [a,b]\nearrow[c,d]\Longrightarrow\DD(a,b)=\DD(c,d),\qquad
 \text{for all }a\leq b\text{ and }c\leq d\text{ in }L.
 \]
For a simple \emph{geometric} lattice (or \emph{combinatorial
geometry}) $L$, the dimension monoid of
$L$ reflects the modularity of $L$, as $\Dim L$ is isomorphic to
$\ZZ^+$ if $L$ is modular, and to $\two$, the two-element semilattice,
otherwise, see F. Wehrung \cite[Corollary~7.12]{WDim}.

A lattice-theoretical antithesis of the topic of continuous
geometries or combinatorial geometries is provided
by \emph{convex geometries}, see K.V. Adaricheva, V.A. Gorbunov, and
V.I. Tumanov \cite{AGT} for a survey. The corresponding algebraic
antithesis of modularity is the \emph{\jsdy}, which is the
quasi-identity
 \begin{equation}\label{Eq:JSD}
 x\vee y=x\vee z\Longrightarrow x\vee y=x\vee(y\wedge z).
 \end{equation}
For a \jsd\ lattice $L$, the substitute of perspectivity for modular
lattices is the relation of \emph{join-dependency} $\DR_L$,
introducing a \emph{polarization} among dimensions, for example,
 \[
 p\DR_Lq\Longrightarrow\DD(p_*,p)+\DD(q_*,q)=\DD(q_*,q),
 \]
for all completely \jirr\ elements $p$ and $q$ of $L$, see
Corollary~\ref{C:tr2ll}. For a general finite lattice~$L$, it is
well-known that if $\tr_L$ denotes the transitive closure of $\DR_L$,
then $\Con L$ is isomorphic to the lattice of lower subsets of the
quasi-ordered system $(\J(L),\tr_L)$,
see R. Freese, J. Je\v{z}ek, and J.B. Nation \cite[Theorem~2.35]{FJN}.
Our methods of computation of the dimension monoid establish the
relation
 \[
 \Dim L\cong\EE(P,\tr_L'),
 \]
the \emph{primitive monoid} generated by the quasi-ordered system
$(P,\tr'_L)$, where $\tr'_L$ is a transitive binary relation on a set
$P$ of smaller size than $\J(L)$, see Theorem~\ref{T:DimBCF}.
If $L$ is \jsd, then it turns out that $P=\J(L)$, furthermore,
$\tr'_L$ and $\tr_L$ are identical (see Corollary~\ref{C:DimBCF}). In
that case the dimensionality on $L$ can be defined by a family of
$\ZZ^+\cup\set{\infty}$-valued functions that can be described
explicitly, see Lemma~\ref{L:Valdp}. These results are also extended
to many infinite lattices, see, for example,
Theorem~\ref{T:DepDimFct}.

As an immediate corollary of our results, we mention the following:
\emph{a finite \jsd\ lattice $L$ is a lower bounded homomorphic image
of a free lattice if{f} $\Dim L$ satisfies the quasi-identity
$2x=x\Rightarrow x=0$}, see Corollary~\ref{C:finLB}.

We finally use these results to extend to the dimension monoid some
known results about the congruence lattice of the \emph{tensor
product}
$A\otimes B$ of lattices $A$ and $B$, see G. Gr\"atzer, H. Lakser, and
R.W. Quackenbush \cite{GLQu81} and G. Gr\"atzer and F. Wehrung
\cite{GrWe2,GrWe3,GrWe4}. For finite lattices
$A$ and $B$, it is not always the case that $\Dim(A\otimes B)$ is
isomorphic to $\Dim A\otimes\Dim B$, however, we prove that related
positive statements hold. For example, we obtain that
\emph{for finite, \jsd\ lattices $A$ and $B$, the relation
$\Dim(A\bp B)\cong\Dim A\otimes\Dim B$ holds} (see
Corollary~\ref{C:AbpB}), where $A\bp B$, the
\emph{box product} of $A$ and $B$, is a variant of the tensor
product, see G. Gr\"atzer and F. Wehrung \cite{GrWe4}. In particular,
we prove that $A\bp B$ is \jsd\ whenever both $A$ and $B$ are \jsd.
This result does not extend to the classical tensor product, see the
counterexample of Section~\ref{S:TsCtEx}.

We conclude the paper with a list of open problems.

\section{Basic notions}\label{S:Basic}

For a partially ordered set $P$ and a subset $X$ of $P$, we put 
 \[
 \dnw X=\setm{p\in P}{\exists x\in X\text{ such that }p\leq x},\ 
 \upw X=\setm{p\in P}{\exists x\in X\text{ such that }x\leq p}.
 \]
For elements $x$ and $y$ of $P$, we write $x\prec y$, if $x<y$ and
no $z\in P$ satisfies that $x<z<y$. We write $x\Vert y$, if
$x\nleq y$ and $y\nleq x$.

For elements $a\leq b$ and $c\leq d$ of a lattice $L$, we write
$[a,b]\nearrow[c,d]$, if $a=b\wedge c$ and $d=b\vee c$. 

We put $L^-=L\setminus\set{0}$ if $L$ has a zero element and
$L^-=L$ otherwise, and we denote by $\J(L)$ (resp., $\JC(L)$) the
set of all \jirr\ (resp., completely \jirr) elements of $L$ (so
$\JC(L)\subseteq\J(L)\subseteq L^-$). For $p\in\JC(L)$, we denote by
$p_*$ the unique lower cover of $p$.

A lattice $L$ is \emph{\jsd} (see \cite{AGT}), if it satisfies the
quasi-identity \eqref{Eq:JSD}. An important class of \jsd\ lattices is
the class of so-called \emph{lower bounded homomorphic
images of free lattices}, that are the
images of finitely generated free lattices under lower bounded lattice
homomorphisms, see \cite{FJN}. 

Let $M$ be a \cm. We say that $M$ is
\begin{itemize}
\item \emph{cancellative}, if $a+c=b+c$ implies that $a=b$, for
all $a$, $b$, $c\in M$;

\item \emph{separative}, if $2a=a+b=2b$ implies that $a=b$, for
all $a$, $b\in M$;

\item \emph{strongly separative}, if $a+b=2b$ implies that $a=b$,
for all $a$, $b\in M$.
\end{itemize}

The notions of cancellativity, separativity, and strong
separativity have even more precise analogues in the
theory of positively preordered monoids. Namely, a positively
preordered monoid is cancellative if{f} it is isomorphic to
the positive cone~$G^+$ of some partially preordered Abelian
group $G$; it is separative if{f} it embeds into a product of
structures of the from $G^+\cup\set{+\infty}$ (see
F. Wehrung \cite{Wehr94}), and strongly separative if{f} it embeds into
a product of structures of the form $G^+\cup\mathcal{R}(I)$, where
$\mathcal{R}(I)$ is the lexicographical sum of $I$ copies of the
real line along a chain $I$ (see C. Moreira dos Santos \cite{MrDS}).

We say that $M$ is a \emph{refinement monoid}, if for all $a_0$,
$a_1$, $b_0$, $b_1\in M$ such that $a_0+a_1=b_0+b_1$, there are
$c_{i,j}$ ($i$, $j<2$) in $M$ such that $a_i=c_{i,0}+c_{i,1}$ and
$b_i=c_{0,i}+c_{1,i}$ for all $i<2$.

The \emph{algebraic preordering} and the \emph{absorption} on $M$ are
respectively defined by
 \begin{align*}
 x\leq y&\Longleftrightarrow\exists z\in M
 \text{ such that }x+z=y,\\
 x\ll y&\Longleftrightarrow x+y=y,
 \end{align*}
for all $x$, $y\in M$. We put $\ZZ^+=\set{0,1,2,\ldots}$ and
$\ZZb=\ZZ^+\cup\set{\infty}$, endowed with its canonical
structure of commutative monoid. We observe that $\ZZb$ is
separative, although not strongly separative.

\section{Dimension functions on a lattice}\label{S:DimFct} 

For a partially ordered set $P$, we put
 \[
 \diag{P}=\setm{(x,y)\in P\times P}{x\leq y}.
 \]

\begin{definition}\label{D:DimFct}
Let $L$ be a lattice, let $M$
be a \cm. A \emph{$M$-valued dimension function} on $L$ is a map
$f\colon\diag{L}\to M$ that satisfies the following equalities, for
all $x$, $y$, $z\in L$:
\begin{itemize}
\item[(D0)] $f(x,x)=0$;

\item[(D1)] $f(x,z)=f(x,y)+f(y,z)$ if $x\leq y\leq z$; 

\item[(D2)] $f(x\wedge y,x)=f(y,x\vee y)$.
\end{itemize}

We let $\R(f)$ denote the submonoid of $M$ generated by the range
of $f$. 

A dimension function $f\colon\diag{L}\to M$ is \emph{universal}, if
for every \cm\ $N$ and every $N$-valued dimension function
$g\colon\diag{L}\to N$, there exists a unique monoid homomorphism
$\varphi\colon M\to N$ such that $g=\varphi\circ f$. We say that $f$
is \emph{separating}, if the restriction of $f$ from $\diag{L}$ to
$\R(f)$ is universal.
\end{definition}

Of course, there is, up to isomorphism, a unique universal
dimension function $f\colon\diag{L}\to M$. The
monoid $M$ is called in \cite{WDim} the \emph{dimension monoid} of
$L$, and denoted by $\Dim L$, while the corresponding dimension
function is denoted by $\DD\colon\diag{L}\to\Dim L$. Hence, the
universal dimension functions on~$M$ are exactly the compositions
with $\DD$ of any monoid embedding from $\Dim L$ into some \cm. 

We shall now present another way to obtain dimension functions on
a lattice~$L$. We shall first present standard definitions
concerning the relation of join-dependency on $L$, see \cite{FJN}.
For any $a\in L$, we put $\J_L(a)=\setm{p\in\J(L)}{p\leq a}$. The
\emph{join-dependency relation} $\DR_L$ on $L\times\J(L)$ is
defined by the rule
 \[
 a\DR_L q\text{ if there exists }x\in L\text{ such that }
 a\leq q\vee x \text{ while }a\nleq y\vee x\text{ for all }y<q
 \]
for all $a\in L$ and all $q\in\J(L)$, and we let $\tr_L$ (resp.,
$\utr_L$) denote the transitive closure (resp., the reflexive and
transitive closure) of $\DR_L$ on $\J(L)$. 
We shall use the notations
$\J(a)$, $\DR$, $\tr$, $\utr$ in case the lattice $L$ is
understood from the context. 

For a binary relation $\alpha$ on
$\J(L)$ and $p\in\J(L)$, we put
 \[
 [p]^\alpha=\setm{q\in\J(L)}{p\mathbin{\alpha}q},\text{ for all
 }p\in\J(L).
 \]
This notation will be used for $\alpha$ being
either $\tr_L$ or $\utr_L$. 

\begin{definition}\label{D:CanDimFct}
Let $L$ be a lattice. For $p\in\J(L)$, we define a map
$d_p\colon\diag{L}\to\nobreak\ZZb$ by the rule 
 \[
 d_p(x,y)=\begin{cases}
 0,&\text{if }\J(x)\cap[p]^{\utr}=\J(y)\cap[p]^{\utr},\\
 1,&\text{if }p\in\J(y)\setminus\J(x)\text{ and }
 \J(x)\cap[p]^{\tr}=\J(y)\cap[p]^{\tr},\\
 \infty,&\text{if }\J(x)\cap[p]^{\tr}\subset\J(y)\cap[p]^{\tr},
 \end{cases}
 \]
for all $(x,y)\in\diag{L}$.
\end{definition}

An immediate computation yields the following special values of
the function $d_p$:

\begin{lemma}\label{L:Valdp}
Let $p$ be a completely \jirr\
element of a lattice $L$. Then $d_p(p_*,p)$ can be computed as
follows:
 \[
 d_p(p_*,p)=\begin{cases}
 1&(\text{if }p\not\tr p),\\
 \infty&(\text{if }p\tr p).
 \end{cases}
 \]
\end{lemma}

We now recall some standard terminology about join-covers, see
again \cite{FJN}. For subsets $X$ and $Y$ of a lattice $L$, we say
that $X$ \emph{refines} $Y$, in notation
$X\sqsubseteq Y$, if for all $x\in X$ there exists $y\in Y$ such
that $x\leq y$ (we do not use the symbol $\ll$, because
here it denotes absorption in \cm s, see Section~\ref{S:Basic}). For
$a\in L^-$, a \emph{join-cover} of $a$ is a nonempty finite subset $X$
of~$L^-$ such that
$a\leq\bigvee X$. We say that a join-cover $X$ of $a$ is
\emph{nontrivial}, if $a\notin\dnw X$. A nontrivial join-cover $X$
of $a$ is \emph{minimal}, if every nontrivial join-cover $Y$ of
$a$ such that $Y\sqsubseteq X$ contains $X$. Observe that $X$ is then
a subset of $\J(L)$.

\begin{definition}\label{D:WJCRP}
We say that a lattice $L$ has
the \emph{weak minimal join-cover refinement property}, if for
every $a\in L^-$, every nontrivial join-cover of $a$ can be
refined to a minimal nontrivial join-cover of $a$.
\end{definition}

The classical definition of the minimal join-cover refinement
property (see \cite{FJN}) is obtained, from the definition above, by
adding the condition that every element has only finitely many
nontrivial join-covers. Hence all finite, or, more generally, finitely
presented lattices, lower bounded homomorphic
images of free lattices, and projective lattices
have the weak minimal join-cover refinement property.

\begin{example}\label{Ex:wmrnotmr}
Let $L$ be the set of all finite subsets $X$ of $\omega$ such that
$\set{1,n}\subseteq X$ implies that $0\in X$, for all $n\geq2$. Then
$L$ is a locally finite, atomistic, \jsd\ lattice with zero.
Furthermore, every finite sublattice of $L$ is a lower bounded homomorphic
image of a free lattice. Since
every principal ideal of $L$ is finite, $L$ has the weak minimal
join-cover refinement property.

However, $L$ does not have the minimal join-cover refinement
property. Indeed, if we put $a=\set{0}$, $b=\set{1}$, and
$b_n=\set{n+2}$ for all $n<\omega$ (these are all the atoms of $L$),
then $\set{b,b_n}$ is a minimal nontrivial join-cover of $a$, and there
are infinitely many such.
\end{example}

For a partially ordered set $P$, an \emph{antichain} of $P$ is a subset
of $P$ whose elements are mutually incomparable. The following lemma
is well-known, see, for example, \cite[Theorem~VIII.2.2]{Birk93}:

\begin{lemma}\label{L:WfantWF}
Let $P$ be a well-founded partially
ordered set. Then the set $\Ant P$ of all finite antichains of
$P$, ordered by $\sqsubseteq$, is well-founded.
\end{lemma}

Now the result mentioned above:

\begin{proposition}\label{P:WF2WJR}
Every well-founded lattice
has the weak minimal join-cover refinement property.
\end{proposition}

\begin{proof}
Let $a\in L^-$, let $X$ be a nontrivial join-cover
of $a$, we prove that there exists a minimal nontrivial join-cover
$Y$ of $a$ such that $Y\sqsubseteq X$. First, by replacing $X$ by
its set $\Max X$ of maximal elements, we may assume without loss
of generality that $X$ is an antichain of $L$. Furthermore, it
follows from Lemma~\ref{L:WfantWF} that $\Ant L$ is well-founded
under $\sqsubseteq$, thus there exists a $\sqsubseteq$-minimal
finite antichain
$Y$ of $L$ such that $a\leq\bigvee Y$ and $Y\sqsubseteq X$. Now
let $Z\sqsubseteq Y$ be a join-cover of $a$. Then $\Max Z$ is a
nontrivial join-cover of $a$ and
$\Max Z\sqsubseteq Y$, thus, by the definition of $Y$, $\Max Z=Y$,
whence $Y\subseteq Z$, so $Y$ is a minimal nontrivial join-cover
of $a$.
\end{proof}

We observe that the lattice of Example~\ref{Ex:wmrnotmr} is
well-founded, although it does not have the minimal join-cover
refinement property.

Now we can state the following result:

\begin{proposition}\label{P:dpDimFct}
Let $L$ be a lattice that
satisfies the weak minimal join-cover refinement property. Then
for all $p\in\J(L)$, the map $d_p$ is a $\ZZb$-valued dimension
function on $L$.
\end{proposition}

\begin{proof}
The items (D0) and (D1) of Definition~\ref{D:DimFct}
are trivially satisfied by $d_p$. Let $x$, $y\in L$ and
$n\in\ZZb$, we prove that $n=d_p(x\wedge y,x)$ if{f}
$n=d_p(y,x\vee y)$. We separate cases: 

\case{1} $n=\infty$. It suffices to prove that
$\J(x)\cap[p]^{\tr}\subseteq\J(y)$ if{f}
$\J(x\vee y)\cap[p]^{\tr}\subseteq\J(y)$. The implication from right to
left is trivial. Conversely, suppose that
$\J(x)\cap[p]^{\tr}\subseteq\J(y)$, and let $q\in\J(x\vee
y)\cap[p]^{\tr}$. Suppose that $q\nleq y$. Hence, by assumption,
$q\nleq x$. Thus, since $q\leq x\vee y$ and by the weak minimal
join-cover refinement property, there exists a nontrivial
join-cover $Z$ of $q$ such that $Z\sqsubseteq\set{x,y}$. For all
$z\in Z$, the relation $q\DR z$ holds, thus $p\tr z$, whence, if
$z\leq x$, we obtain that $z\in\J(x)\cap[p]^{\tr}\subseteq\J(y)$, thus
$z\leq y$. Therefore, $Z=(Z\cap\dnw x)\cup(Z\cap\dnw y)$ is
contained in $\dnw y$, whence $q\leq y$, a contradiction.

So, $q\leq y$, which completes the proof of the assertion of
Case~1. 

\case{2} $n<\infty$. By the result of Case~1, it suffices to prove
that if $\J(x)\cap[p]^{\tr}\subseteq\J(y)$, then
$p\in\J(x)\setminus\J(x\wedge y)$ if{f} $p\in\J(x\vee
y)\setminus\J(y)$. The implication from left to right is trivial.
Conversely, suppose that $p\in\J(x\vee y)\setminus\J(y)$. So
$p\notin\J(x\wedge y)$; suppose that $p\notin\J(x)$. Let $Z$ be a
minimal nontrivial join-cover of $p$ such that
$Z\sqsubseteq\set{x,y}$. For all $z\in Z$, the relation
$p\DR z$ holds, thus, if $z\in Z\cap\dnw x$, then
$z\in\J(x)\cap[p]^{\tr}\subseteq\J(y)$. Therefore, $Z=(Z\cap\dnw
x)\cup(Z\cap\dnw y)$ is contained in $\dnw y$, whence
$p\leq\bigvee Z\leq y$, a contradiction. So $p\leq x$, which
completes the proof of the assertion of Case~2.
\end{proof}

\section{Join-dependency and dimension in BCF lattices}\label{S:JDD}

We shall use the following terminology, introduced in \cite{WDim}.

\begin{definition}\label{D:BCF}
A partially ordered set $P$ is \emph{BCF}, if every bounded chain of
$P$ is finite.
\end{definition}

Let $L$ be a lattice. We introduce on $\JC(L)$ the following
refinements of the join-dependency relation $\DR=\DR_L$ on $L$.

\begin{notation}\label{No:VariousD}
For any $p$, $q\in\JC(L)$ such
that $p\neq q$, we write that \begin{itemize}
\item $p\DR^0q$, if there exists $x\in L$ such that $p\nleq x$,
$p_*\vee q_*\leq x$, and $p\vee x=q\vee x$. 

\item $p\DR^1q$, if there exists $x\in L$ such that $p\nleq x$,
$q_*\leq x$, $p\leq q\vee x$, and either $p\vee x=p_*\vee x$ or
$q\nleq p\vee x$. 

\item $p\DR^\infty q$, if there exists $x\in L$ such that $p\nleq
x$, $p\vee x=q\vee x$, $p_*\vee q_*\leq x$, and
$p\nleq q\vee(x\wedge(p\vee q))$.
\end{itemize}
The various possibilities
are illustrated on Figure~1.
\end{notation}

\begin{figure}[htb]
\includegraphics{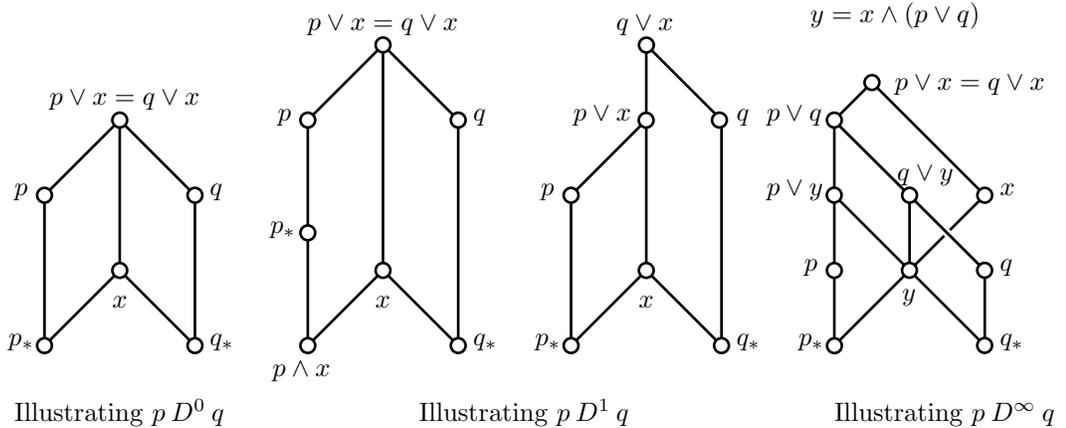}
\caption{The relations $\DR^0$, $\DR^1$, and $\DR^\infty$}
\end{figure}

The following lemma expresses the main relations between $\DR^0$,
$\DR^1$, $\DR^\infty$, the join-dependency relation $\DR$, and
dimension. For all $p\in\JC(L)$, we put $\DD(p)=\DD(p_*,p)$. 

\begin{lemma}\label{L:VariousD}
For all distinct $p$, $q\in\JC(L)$,
the following assertions hold:
\begin{enumerate}
\item $p\DR q$ if{f} either $p\DR^0q$ or $p\DR^1q$; 

\item $p\DR^\infty q$ implies that $p\DR^0q$; 

\item $p\DR^0q$ implies that $\DD(p)=\DD(q)$; 

\item $p\DR^1q$ implies that $\DD(p)\ll\DD(q)$; 

\item $p\DR^\infty q$ implies that $\DD(p)=\DD(q)=2\DD(p)$; 

\item If $L$ is \jsd, then $p\DR^0q$ never holds.
\end{enumerate}

\end{lemma}

\begin{proof} (i) If $x$ is a witness for either $p\DR^0q$ or
$p\DR^1q$, then $p\leq q\vee x$ and $p\nleq x=q_*\vee x$, whence
$p\DR q$. Conversely, suppose that $p\DR q$, thus there exists
$x\geq q_*$ in $L$ such that $p\leq q\vee x$ and $p\nleq x$. If
$p_*\vee x=p\vee x$, then $x$ witnesses that $p\DR^1q$. Suppose
now that $p\nleq p_*\vee x$. Then we may replace $x$ by $p_*\vee
x$, and thus suppose that $p_*\vee q_*\leq x$. If $p\vee x=q\vee
x$, then $x$ witnesses that $p\DR^0q$. If $p\vee x<q\vee x$, then
$x$ witnesses that $p\DR^1q$.

(ii) is obvious.

(iii) {}From $p\nleq x$ and $p_*\leq x$ follows that $p\wedge
x=p_*$, thus $[p_*,p]\nearrow[x,p\vee x]$, whence
$\DD(p)=\DD(x,p\vee x)$. Similarly, $\DD(q)=\DD(x,q\vee x)$. The
conclusion follows from $p\vee x=q\vee x$. 

(iv) Suppose first that $p\vee x=p_*\vee x$. Hence the elements
$p>p_*$ and $x$ generate a pentagon with bottom $p\wedge x$ and
top $p\vee x$, which yields the following relation:
 \begin{equation}\label{Eq:DDppveex}
 \DD(p)\ll\DD(x,p\vee x).
 \end{equation}
Moreover, $q\nleq x$ (otherwise $p\leq q\vee x=x$,
a contradiction) and $q_*\leq x$, thus $[q_*,q]\nearrow[x,q\vee
x]$, so we obtain the equality 
 \begin{equation}\label{Eq:DDqqveex}
 \DD(q)=\DD(x,q\vee x).
 \end{equation}
{}From $x\leq p\vee x\leq q\vee x$ follows that
$\DD(x,p\vee x)\leq\DD(x,q\vee x)$, therefore, by
\eqref{Eq:DDppveex} and \eqref{Eq:DDqqveex}, we obtain the desired
relation $\DD(p)\ll\DD(q)$. 

Suppose now that $p\nleq p_*\vee x$ and $q\nleq p\vee x$. After
replacing $x$ by $p_*\vee x$, we may assume that $p_*\leq x$, so
$[p_*,p]\nearrow[x,p\vee x]$, which yields the following
equality: 
 \begin{equation}\label{Eq:DDp=pveex}
 \DD(p)=\DD(x,p\vee x).
 \end{equation}
Moreover, since $q\nleq p\vee x$, the elements
$p\vee x>x$ and $q$ generate a pentagon with bottom $q_*$ and top
$q\vee x$, which yields the following relation:
 \begin{equation}\label{Eq:DDxpvxllq}
 \DD(x,p\vee x)\ll\DD(q).
 \end{equation}
{}From \eqref{Eq:DDp=pveex} and
\eqref{Eq:DDxpvxllq} follows again that $\DD(p)\ll\DD(q)$.

(v) Let $x$ witness that $p\DR^\infty q$. We put $y=x\wedge(p\vee
q)$. It follows from item (ii) above that $p\DR^0q$, whence, by
(iii), $\DD(p)=\DD(q)$. Furthermore, $[q_*,q]\nearrow[x,q\vee x]$
and $[q_*,q]\nearrow[y,q\vee y]$, which yields the equalities 
 \begin{equation}\label{Eq:DDqxqp}
 \DD(q)=\DD(x,q\vee x)=\DD(y,q\vee y).
 \end{equation}
Furthermore, from the relations
$[y,p\vee q]\nearrow[x,q\vee x]$ and $[p_*,p]\nearrow[q\vee y,p\vee q]$
follows that 
 \begin{equation}\label{Eq:bunchofDs}
 \DD(x,q\vee x)=\DD(y,p\vee q)=\DD(y,q\vee y)+\DD(q\vee y,p\vee q)
 =\DD(q)+\DD(p).
 \end{equation}
The desired conclusion follows from
\eqref{Eq:DDqxqp} and \eqref{Eq:bunchofDs}.

(vi) Suppose that $L$ is \jsd\ and that $x$ witnesses that
$p\DR^0q$, \emph{i.e.}, $p\nleq x$, $p_*\vee q_*\leq x$, and
$p\vee x=q\vee x$. It follows from the \jsdy\ of $L$ that $p\vee
x=(p\wedge q)\vee x$, but $p\nleq q$ (otherwise, since $p\neq q$,
$p\leq q_*\leq x$, a contradiction), thus $p\wedge q\leq p_*$,
whence $(p\wedge q)\vee x=x$ (because $p_*\leq x$), so $p\leq
p\vee x=(p\wedge q)\vee x=x$, a contradiction.
\end{proof}

\begin{corollary}\label{C:tr2ll}
Let $L$ be a lattice, let $p$, $q\in\JC(L)$. We denote by $\trc$ the
transitive closure of the restriction of the join-dependency relation
to $\JC(L)$. Then the following assertions hold:
\begin{enumerate}
\item If $L$ has the weak minimal join-cover refinement property,
then $\DD(p)\ll\nobreak\DD(q)$ implies that $p\tr q$.

\item If $L$ is \jsd, then $p\trc q$ implies that $\DD(p)\ll\DD(q)$.
\end{enumerate}
\end{corollary}

\begin{proof}
(i) It follows from Proposition~\ref{P:dpDimFct} that
$d_p(p_*,p)\ll d_p(q_*,q)$, but $d_p(p_*,p)>\nobreak0$ (see
Lemma~\ref{L:Valdp}), whence
$d_p(q_*,q)=\infty$, \emph{i.e.}, by the definition of $d_p$,
$\J(q_*)\cap[p]^{\tr}\subset\J(q)\cap[p]^{\tr}$. Thus
$q\in[p]^{\tr}$, \emph{i.e.}, $p\tr q$.

(ii) It suffices to consider
the case where $p\DR q$, in which case, by
Lemma~\ref{L:VariousD}(i,vi), $p\DR^1q$, whence, by
Lemma~\ref{L:VariousD}(iv), $\DD(p)\ll\DD(q)$.
\end{proof}

For a finite lattice $L$, the congruence lattice $\Con L$ of $L$
can be computed from the $\utr$ relation on $\J(L)$, see
\cite[Theorem~2.35]{FJN}. 

For a BCF lattice $L$ with zero, our following result gives a related
way to compute $\Dim L$.

\begin{theorem}\label{T:DimBCF}
Let $L$ be a BCF lattice with zero. Then
$\Dim L$ is isomorphic to the \cm\ $\Dim'L$ defined by generators
$\DD'(p)$, for $p\in\J(L)$, and the following relations:
 \begin{align*}
 \DD'(p)&=\DD'(q)&&\text{if }p\DR^0q,\\
 \DD'(p)&\ll\DD'(q)&&\text{if either }p\DR^1q
 \text{ or }p\DR^\infty q,
 \end{align*}
for all $p$, $q\in\J(L)$. The
isomorphism carries $\DD(p)$ to $\DD'(p)$, for all $p\in\J(L)$.
\end{theorem}

In particular, the assumption of Theorem~\ref{T:DimBCF} holds for
any \emph{finite} lattice $L$.

\begin{proof}
We first observe that $L$ is well-founded, thus, by
Lemma~\ref{P:WF2WJR}, $L$ has the weak minimal join-cover refinement
property (this is easily seen to fail as a rule for BCF lattices
without zero). Furthermore, since every bounded interval of
$L$ is n{\oe}therian (\emph{i.e.}, dually well-founded), every \jirr\
element of $L$ is completely \jirr\ (\emph{i.e.}, $\J(L)=\JC(L)$).

It follows from Lemma~\ref{L:VariousD} that there
exists a monoid homomorphism\linebreak
$\varphi\colon\Dim'L\to\Dim L$ such that
$\varphi(\DD'(p))=\DD(p)$, for all $p\in\J(L)$. To prove the
converse, we use the alternate presentation of $\Dim L$ \emph{via}
caustic pairs given in \cite[Chapter~7]{WDim}. More precisely,
$\Dim L$ is defined by generators $\DD(a,b)$, where $a\prec b$
in~$L$, subjected to relations given by (7.1)--(7.3) of
\cite[page 318]{WDim}. 

For all $x\prec y$ in $L$, there exists, since $L$ is
well-founded, a minimal element $p\in L$ such that $p\leq y$ and
$p\nleq x$. {}From the minimality assertion on $p$ follows that
$p\in\J(L)$. The assumption that $p_*\nleq x$ would contradict the
minimality assumption on $p$, thus $p_*\leq x$ and then $p\wedge
x=p_*$. Moreover, from $x<p\vee x\leq y$ and $x\prec y$ follows
that $p\vee x=y$, so, finally, $[p_*,p]\nearrow[x,y]$. 

Furthermore, if $q\in\J(L)$ such that $[q_*,q]\nearrow[x,y]$,
then $p\DR^0q$, whence $\DD'(p)=\DD'(q)$. This entitles us to
define, for all $x\prec y$ in
$L$, $\DD'(x,y)=\DD'(p)$ for any $p\in\J(L)$ such that
$[p_*,p]\nearrow[x,y]$. We shall prove that the map $\DD'$ thus
defined on all pairs $(x,y)\in L\times L$ such that $x\prec y$
satisfies the equations listed in (7.1)--(7.3) of
\cite[page 318]{WDim}. It is convenient to start with the following
easy claim. 

\begin{sclaim}
Let $a\prec b$ and $c\prec d$ in $L$. If
$[a,b]\nearrow[c,d]$, then $\DD'(a,b)=\DD'(c,d)$.
\end{sclaim}

\begin{scproof}
Let $p\in\J(L)$ such that $[p_*,p]\nearrow[a,b]$.
Then $[p_*,p]\nearrow[c,d]$ as well, whence
$\DD'(a,b)=\DD'(c,d)=\DD'(p)$.
\end{scproof}

To verify the relations (7.1)--(7.3) of \cite[page 318]{WDim}
amounts to verifying the following cases.
\smallskip

\noindent\textbf{The relations (7.1).} We are given elements $u$,
$v$, $x$, and $y$ of $L$ such that $u\prec\nobreak x<\nobreak v$,
$u<y\prec v$,
$x\wedge y=u$, and $x\vee y=v$. We need to verify that
$\DD'(u,x)=\DD'(y,v)$. This is obvious by the claim above since
$[u,x]\nearrow[y,v]$.
\smallskip

\noindent\textbf{The relations (7.2).} We are given elements $u$,
$v$, $x$, $y$, $z$, and $t$ of $L$ such that
$u\prec x\leq y\prec z<v$, $u<t\prec v$, $t\wedge z=u$, and
$t\vee x=v$. We need to verify that $\DD'(y,z)\ll\DD'(u,x)$. So, let
$p$, $q\in\J(L)$ such that $[p_*,p]\nearrow[y,z]$ and
$[q_*,q]\nearrow[u,x]$, we need to prove that
$\DD'(p)\ll\DD'(q)$. 

Suppose first that $p_*\nleq u$. If $p_*\leq t$, then $p_*\leq
t\wedge z=u$, a contradiction; whence $p_*\nleq t$. But $t\prec v$
and $p\leq z<v$, whence $p_*\vee t=v$. Moreover, $p\nleq t$
(otherwise $p\leq t\wedge z=u$, a contradiction), $q\nleq t$
(otherwise $q\leq t\wedge x=u$, a contradiction), thus $p\vee
t=q\vee t=v$. Hence $t$ witnesses that $p\DR^1q$, whence
$\DD'(p)\ll\DD'(q)$. 

Now suppose that $p_*\leq u$. Then $p\nleq t$ (otherwise $p\leq
t\wedge z=u\leq y$, a contradiction), $p\vee t=q\vee t=v$,
$p_*\vee q_*\leq u\leq t$. Furthermore,
$t\wedge(p\vee q)\leq t\wedge z\leq u$ and $p\nleq x=q\vee u$,
whence $p\nleq q\vee(t\wedge(p\vee q))$. Therefore, $t$ witnesses
that $p\DR^\infty q$, so, again, $\DD'(p)\ll\DD'(q)$. \smallskip

\noindent\textbf{The relations (7.3).} We are given elements $u$,
$v$, $x$, $y$, $z$, and $t$ of $L$ such that $u<z\prec y\leq
x\prec v$, $u\prec t<v$, $x\wedge t=u$, $z\vee t=v$. We need to
verify that $\DD'(z,y)\ll\DD'(x,v)$. Let $p$, $q\in\J(L)$ such
that $[p_*,p]\nearrow[z,y]$ and $[q_*,q]\nearrow[u,t]$, we need
to verify that $\DD'(p)\ll\DD'(q)$. We observe that $q\vee z=q\vee
u\vee z=t\vee z=v$, whence $p\leq q\vee z$. If $q\leq p\vee z$,
then $q\leq y$, but
$q\wedge y=q\wedge t\wedge y=q\wedge u\leq u$, a contradiction;
hence $q\nleq p\vee z$. Furthermore, $q_*\leq u\leq z$ and $p\nleq
z$. Therefore, $z$ witnesses that $p\DR^1q$, so, again,
$\DD'(p)\ll\DD'(q)$.
\end{proof}

As an immediate consequence of Theorem~\ref{T:DimBCF} and
Lemma~\ref{L:VariousD}, we obtain the following: 

\begin{corollary}\label{C:DimBCF}
Let $L$ be a BCF \jsd\ lattice with zero.
Then $\Dim L$ is the \cm\ defined by generators $\olp$ (for
$p\in\J(L)$) and relations $\olp\ll\olq$ for all $p$, $q\in\J(L)$
such that $p\DR q$ (resp., $p\tr q$).
\end{corollary}

We observe that if a \jsd\ lattice $L$ is BCF, then every interval of
$L$ is \emph{finite}, see, for example, \cite[Proposition~3.2]{AGT}
or \cite[Theorem~5.59]{FJN}.


For lattices with zero, the result of Corollary~\ref{C:DimBCF} is
weaker than the result of Theorem~\ref{T:DepDimFct}.

\begin{remark}\label{Rk:DinftyNec}
By identifying both $p_*$ and
$q_*$ with zero in the diagram illustrating~$\DR^\infty$ in
Figure~1, we obtain a nine element lattice with \jirr\ elements
$p$ and~$q$ with $p\tr q$ (and even $p\DR^0q$) although
$\DD(p)\not\ll\DD(q)$. Hence the use of $\DR^\infty$ is necessary
in the statement of Theorem~\ref{T:DimBCF}.
\end{remark}

\begin{example}\label{Ex:CoP}
For a partially ordered set $P$, we denote by $\Co(P)$ the
lattice of all subsets $X$ of $P$ that are \emph{order-convex},
\emph{i.e.}, $u\leq p\leq v$ and $\set{u,v}\subseteq X$ implies that
$p\in X$, for all $u$, $v$, $p\in P$. The lattices of the form
$\Co(P)$ are studied in G. Birkhoff and M.K. Bennett \cite{BB}, where
it is proved, in particular, that $\Co(P)$ is \jsd. Of course, the
completely \jirr\ elements of $\Co(P)$ are the singletons of elements
of $P$, and the $\DR$ relation on these is given by
$\set{p}\DR\set{q}$ if{f} $p\dr q$, where $\dr$ is the binary relation
on $P$ given by the rule
 \[
 p\dr q\qquad\text{if{f}}\qquad
 \exists r\in P\text{ such that either }q<p<r\text{ or }r<p<q,
 \quad\text{for all }p,\,q\in P.
 \]
It follows from Corollary~\ref{C:DimBCF} that for finite $P$, the
dimension monoid of $\Co(P)$ is the commutative monoid defined by
the generators $\ol{p}$, for $p\in P$, and the relations
$\ol{p}\ll\ol{q}$ for $p$, $q\in P$ such that $p\dr q$.
\end{example}

\section{Primitive monoids}\label{S:Prim} 

We refer to R.S. Pierce \cite[Sections 3.4--3.6]{Pier} for basic
information about primitive refinement monoids. For a nonzero element
$p$ in a \cm\ $M$, we say that $p$ is \emph{\pind}, if $p=x+y$ implies
that either $x=p$ or $y=p$, for all $x$, $y\in M$. A refinement monoid
$M$ is \emph{primitive}, if it is generated as a monoid by its set
of \pind\ elements and its algebraic preordering is antisymmetric. 

Primitive monoids can be constructed as follows. We say that a
\emph{QO-system} is a pair $(P,\tr)$, where $\tr$ is a transitive
binary relation on a set $P$. For a QO-system $(P,\tr)$, let
$\EE(P,\tr)$ denote the \cm\ defined by generators $\olp$, for
$p\in P$, subjected to the relations $\olp\ll\olq$ (\emph{i.e.},
$\olp+\olq=\olq$) for all $p$, $q\in P$ such that $p\tr q$. Then
the primitive monoids are exactly the monoids of the form
$\EE(P,\tr)$ for a QO-system $(P,\tr)$, see
\cite[Proposition~3.5.2]{Pier}; in addition, one can take $\tr$
antisymmetric. The \pind\ elements of $\EE(P,\tr)$ are exactly the
elements $\olp$ for $p\in P$. 

We recall two well-known lemmas about primitive monoids: 

\begin{lemma}[see
{\cite[Proposition~3.4.4]{Pier}}]\label{L:NormForm}
Let $M$ be a primitive monoid. Then every element $a\in M$ has a
unique representation
 \[
 a=\sum_{i<n}a_i,
 \]
in which $n<\omega$, the elements $a_0$, \dots, $a_{n-1}$ are
\pind, and $a_i\not\ll a_j$ for all $i$, $j<n$ with $i\neq j$.
\end{lemma}

We shall call the decomposition of $a$ given in
Lemma~\ref{L:NormForm} the \emph{canonical decomposition} of $a$.

For our next lemma, for a QO-system $(P,\tr)$, we denote by
$\FF(P,\tr)$ the set of all mappings $\lx\colon P\to\ZZb$ such
that $p\tr q$ implies that $\lx(q)\ll\lx(p)$, for all $p$,
$q\in P$. Of course, $\FF(P,\tr)$ is an additive submonoid of
$(\ZZb)^P$. For any $p\in P$, we denote by $\otp$ the element of
$\FF(P,\tr)$ defined by the rule 
 \[
 \otp(q)=\begin{cases}
 \infty,&\text{if }q\tr p,\\ 1,&\text{if }q=p\ntr p,\\ 0,&\text{if
 }q\nutr p,
 \end{cases}
 \]
for all $q\in P$. We warn the reader that the notation
$\FF(P,\tr)$ used here does not mean the same as the corresponding
notation in \cite[Chapter~6]{WDim}.

It is clear that $p\tr q$ implies that $\otp\ll\otq$, for all $p$,
$q\in P$. In fact, much more can be said, see
\cite[Proposition~6.8]{WDim}: 

\begin{lemma}\label{L:E(P)subF(P)}
There exists a unique monoid
homomorphism from $\EE(P,\tr)$ to $\FF(P,\tr\nobreak)$ that sends
$\olp$ to $\otp$ for all $p\in P$, and it is a monoid embedding.
\end{lemma}

Hence, from now on we shall identify $\EE(P,\tr)$ with its image
under the natural embedding into $\FF(P,\tr)$, thus we will also
identify $\otp$ with $\olp$. Observe that this way, $\EE(P,\tr)$
becomes a submonoid of a direct power of $\ZZb$.
\emph{In particular, $\EE(P,\tr)$ is separative}. 

\begin{lemma}\label{L:StrSepPrim}
Let $(P,\tr)$ be a QO-system.
Then the following are equivalent:
\begin{enumerate}
\item $\EE(P,\tr)$ is strongly separative; 

\item $\EE(P,\tr)$ satisfies the axiom $(\forall
x)(2x=x\Rightarrow x=0)$; 

\item the binary relation $\tr$ is irreflexive.
\end{enumerate}

\end{lemma}

\begin{proof}
(i)$\Rightarrow$(ii) is obvious.

(ii)$\Rightarrow$(iii) Suppose that the assumption of (ii) is
satisfied. Thus, for all $p\in P$, $2\olp\neq\olp$, whence $p\ntr
p$. 

(iii)$\Rightarrow$(i) Suppose that $\tr$ is irreflexive, in
particular, we may identify $P$ with the set of all \pind\
elements of $\EE(P,\tr)$. Let $\la$, $\lb\in\EE(P,\tr)$ such that
$\la+\lb=2\lb$, we prove that $\la=\lb$. Let
 \[
 \la=\sum_{i<m}\la_i\text{ and }\lb=\sum_{j<n}\lb_j
 \]
be the
canonical decompositions of $\la$ and $\lb$, see
Lemma~\ref{L:NormForm}. We put
 \begin{align*}
 X&=\setm{i<m}{\exists j<n\text{ such that
 }\la_i\ll\lb_j},\\ Y&=\setm{i<m}{\exists j<n\text{ such that
 }\lb_j\ll\la_i},
 \end{align*}
and $Z=m\setminus(X\cup Y)$.
Observe that $\lb=\sum_{j<n}\lb_j$ is the canonical decomposition
of~$\lb$, whence $X\cap Y=\es$; thus $X$, $Y$, and $Z$ are pairwise
disjoint. Let $i_0\in Y$. Then there exists $j_0<n$ such that
$\lb_{j_0}\ll\la_{i_0}$, so
$\la+\lb=\sum_{i<m}\la_i+\sum_{j<n,\,j\neq j_0}\lb_j$, and from
the expression on the right hand side of that equality we can
extract (by removing $\lx$ from $\lx+\ly$ whenever $\lx\ll\ly$) a
canonical decomposition of $\la+\lb$ in which $\lb_{j_0}$ does not
occur. However, $2\lb=\sum_{j<n}(\lb_j+\lb_j)$ is the canonical
decomposition of $2\lb$ and $\lb_{j_0}$ occurs there, a
contradiction. Hence $Y=\es$. 

Furthermore, $\la_i\ll\lb$ for all $i\in X$, whence
$\la+\lb=\sum_{i\in Z}\la_i+\sum_{j<n}\lb_j$ is the canonical
decomposition of $\la+\lb$. Since $2\lb=\sum_{j<n}(\lb_j+\lb_j)$
is the canonical decomposition of $2\lb$, it follows from
Lemma~\ref{L:NormForm} that $\sum_{i\in Z}\la_i=\lb$, whence
$\la=\sum_{i\in Z}\la_i+\sum_{i\in X}\la_i=
\lb+\sum_{i\in X}\la_i=\lb$.
\end{proof}

\section{The dependency dimension function on a
lattice}\label{S:trDimFct} 

Until Corollary~\ref{C:inFP}, we shall fix a lattice $L$ which has
the weak minimal join-cover refinement property, see
Definition~\ref{D:WJCRP}. For all $(x,y)\in\diag{L}$, we put
 \begin{equation}\label{Eq:DefDF}
 \lDD(x,y)=(d_p(x,y))_{p\in\J(L)},
 \end{equation}
see Definition~\ref{D:CanDimFct}. So, $\lDD(x,y)$
is an element of $(\ZZb)^{\J(L)}$. The following lemma says more: 

\begin{lemma}\label{L:inFP}
$\lDD(x,y)$ belongs to $\FF(\J(L),\tr)$, for all
$(x,y)\in\diag{L}$.
\end{lemma}

\begin{proof}
Put $\la=\lDD(x,y)$, it suffices to prove that
$\la(p)<\infty$ implies that $\la(q)=0$, for all $p$, $q\in\J(L)$
such that $p\tr q$. By assumption,
$\J(x)\cap[p]^{\tr}=\J(y)\cap[p]^{\tr}$, whence, since $p\tr q$,
$\J(x)\cap[q]^{\utr}=\J(y)\cap[q]^{\utr}$, so $\la(q)=0$, indeed.
\end{proof}

It is convenient to record as follows the immediate consequence of
Lemma~\ref{L:inFP} and Proposition~\ref{P:dpDimFct}: 

\begin{corollary}\label{C:inFP} The map $\lDD$ is a
$\FF(\J(L),\tr)$-valued dimension function on $L$. Furthermore,
$\lDD(p_*,p)=\otp$ for all $p\in\JC(L)$.
\end{corollary}

We shall call $\lDD$ the \emph{dependency dimension function} on
$L$. 

We shall now investigate conditions under which $\lDD$ is
\emph{separating}, see Section~\ref{S:DimFct}.

\begin{notation}
Let $L$ be a lattice. For $(a,b)\in\diag{L}$, we
write that $a\lessdot b$, if there exists $p\in\JC(L)$ such that
$[p_*,p]\nearrow[a,b]$.
\end{notation}

The following lemma shows that the relation $a\lessdot b$ is not
uncommon: 

\begin{lemma}\label{L:manylessdot}
Let $L$ be a lattice. If $L$ is \emph{spatial}, \emph{i.e.}, every
element of $L$ is a join of completely \jirr\ elements of $L$,
then $a\prec b$ implies that $a\lessdot b$, for all $a$, $b\in L$.
\end{lemma}

We observe that the assumption of $L$ being spatial holds for $L$
\emph{dually algebraic} (see Theorem I.4.22 in G. Gierz
\emph{et al.} \cite{Comp}, or Lemma~1.3.2 in V.A. Gorbunov
\cite{Gorb}), thus, in particular, for $L$ \emph{well-founded}.

\begin{proof}
Let $p$ be a completely \jirr\ element of $L$ such that $p\leq b$ and
$p\nleq a$. If $p\wedge a<p_*$, then $p_*\nleq a$, which contradicts
the minimality assumption on~$p$. Therefore, $[p_*,p]\nearrow[a,b]$.
\end{proof}

Now we are ready to prove the main result of this section: 

\begin{theorem}\label{T:DepDimFct}
Let $L$ be a lattice satisfying
the following properties:
\begin{enumerate}
\item the weak minimal join-cover refinement property; 

\item for all $a<b$ in $L$, there are a positive integer $n$ and a
chain 
 \[
 a=x_0\lessdot x_1\lessdot\cdots\lessdot x_n=b,\quad
 \text{for elements }x_0,\dots,x_n\in L.
 \]
\end{enumerate}
Then $\J(L)=\JC(L)$ and the range of $\lDD$ in $L$ generates
$\EE(\J(L),\tr\nobreak)$. Furthermore, if~$L$ is \jsd, then $\lDD$ is
separating (see Definition~\textup{\ref{D:DimFct}}); hence
$\Dim L\cong\EE(\J(L),\tr)$.
\end{theorem}

\begin{proof}
By applying (ii) to the case where $b\in\J(L)$, we immediately
obtain that $\J(L)=\JC(L)$.

For $a$, $b\in L$ such that $a\lessdot b$ and
$p\in\J(L)$ such that $[p_*,p]\nearrow[a,b]$, it follows from
Corollary~\ref{C:inFP} that the equality $\lDD(a,b)=\lDD(p_*,p)$
holds. In particular, $\lDD(a,b)$ belongs to $\EE(\J(L),\tr)$.
Then it follows immediately from assumption (ii) that $\lDD(a,b)$
belongs to $\EE(\J(L),\tr)$ for all $(a,b)\in\diag{L}$. 

Suppose now that $L$ is \jsd. We put $\lDD(p)=\lDD(p_*,p)$, for
all $p\nobreak\in\nobreak\J(L)$. By the paragraph above, there
exists a monoid homomorphism $\pi\colon\Dim L\to\EE(\J(L),\tr)$
such that $\pi(\DD(p))=\lDD(p)$ for all $p\in\J(L)$. To prove the
converse, it suffices, by using the definition of
$\EE(\J(L),\tr)$ \emph{via} generators and relations, to prove
that $p\tr q$ implies that $\DD(p)\ll\DD(q)$, for all $p$,
$q\in\J(L)$. However, this follows immediately from
Corollary~\ref{C:tr2ll}(ii).
\end{proof}

We observe that the assumptions underlying Theorem~\ref{T:DepDimFct}
are not uncommon, for example, they are obviously satisfied by
the lattice $\CB(E)$ of all convex polytopes of any real affine space
$E$. Then Theorem~\ref{T:DepDimFct} yields immediately that for
nontrivial~$E$, $\Dim\CB(E)$ is the two-element semilattice,
which is also easy to verify directly.

As an immediate consequence of Theorem~\ref{T:DepDimFct},
Proposition~\ref{P:WF2WJR}, and Lemma~\ref{L:manylessdot}, we
obtain the following:

\begin{corollary}\label{C:DepDimFct}
Let $L$ be a \jsd, well-founded
lattice in which for all $a<b$ there are a positive integer $n$
and a chain $a=x_0\prec x_1\prec\cdots\prec x_n=b$. Then $L$
satisfies the conditions of Theorem~\textup{\ref{T:DepDimFct}};
whence $\Dim L\cong\EE(\J(L),\tr)$.
\end{corollary}

We observe that the conditions of Corollary~\ref{C:DepDimFct} are
satisfied for $L$ a BCF lattice with zero.

It is worthwhile to record the following consequence of
Lemma~\ref{L:StrSepPrim} and Theorem~\ref{T:DepDimFct}: 

\begin{corollary}\label{C:infLB}
Let $L$ be a \jsd\ lattice
satisfying the assumptions of Theorem~\textup{\ref{T:DepDimFct}}.
Then the following are equivalent:
\begin{enumerate}
\item The join-dependency relation on $\J(L)$ has no cycles. 

\item $\Dim L$ is strongly separative.

\item $\Dim L$ satisfies the axiom
$(\forall x)(2x=x\Rightarrow x=0)$.
\end{enumerate}
\end{corollary}

%

In particular, for a finite lattice $L$, it is well-known (see
\cite{FJN}) that $L$ has no $\DR_L$-cycles if{f} $L$ is a
\emph{lower bounded homomorphic image of a free lattice}. Hence we
obtain the following dimension-theoretical characterization of lower
boundedness:

\begin{corollary}\label{C:finLB}
Let $L$ be a finite \jsd\
lattice. Then the following are equivalent:
\begin{enumerate}
\item $L$ is a lower bounded homomorphic
image of a free lattice.

\item $\Dim L$ is strongly separative.

\item $\Dim L$ satisfies the axiom
$(\forall x)(2x=x\Rightarrow x=0)$.
\end{enumerate}
\end{corollary}

Corollary~\ref{C:finLB} does not extend to lattices that are not
\jsd. For example, for a finite \emph{modular} lattice $L$, the
dimension monoid $\Dim L$ is always \emph{cancellative} (see
\cite[Proposition~5.5]{WDim}), thus \emph{a fortiori} strongly
separative. However, if $L$ is non-distributive, then, since $L$
is modular, it cannot be a lower bounded homomorphic
image of a free lattice. 

In particular, we obtain a well-known result of A. Day, see
\cite[Theorem~2.64]{FJN}:

\begin{theorem}\label{T:Day}
A finite, lower bounded homomorphic
image of a free lattice is an upper bounded homomorphic
image of a free lattice if{f} it is \msd.
\end{theorem}


\begin{proof}
We prove the nontrivial direction. Let $L$ be a
finite lower bounded homomorphic image of a free lattice. It follows
from Corollary~\ref{C:finLB} that $\Dim L$ is strongly separative.
If, in addition, $L$ is \msd, then, since $L$ and its dual lattice have
isomorphic dimension monoids, it follows again from
Corollary~\ref{C:finLB} that~$L$ is an upper bounded homomorphic
image of a free lattice.
\end{proof}

\section{The canonical map from $\Dim A\otimes\Dim B$ to
$\Dim(A\otimes B)$}\label{S:DimTens}


We recall the definition of the \emph{tensor product} of
lattices $A$ and $B$ with zero, see~\cite{GrWe2}. A subset $I$ of
$A\times B$ is a \emph{bi-ideal}, if it contains
$\bot_{A,B}=(A\times\set{0_B})\cup(\set{0_A}\times B)$, and 
 \[
 ((a,x)\in I\text{ and }(a,y)\in I)\Rightarrow(a,x\vee y)\in I
 \text{ for all }x,\,y\in B,
 \]
and symmetrically. Important examples of bi-ideals are the
following:
\begin{itemize}
\item the \emph{pure tensors},
$a\otimes b=\bot_{A,B}\cup
\setm{(x,y)\in A\times B}{x\leq a\text{ and }y\leq b}$, for
$(a,b)\in A\times B$.

\item the \emph{mixed tensors}, \emph{i.e.}, the subsets of
$A\times B$ of the form $(a\otimes b')\cup(a'\otimes b)$, for
$a\leq a'$ in $A$ and $b\leq b'$ in $B$.
\end{itemize}

We denote by $A\ootimes B$ the set of all bi-ideals of $A\times
B$, partially ordered under containment. So $A\ootimes B$ is an
algebraic lattice, we denote by $A\otimes B$ its \jzs\ of compact
elements. 

For lattices $A$ and $B$ with zero, $A\otimes B$ is not always a
lattice, even for $A$ finite, see \cite{GrWe3}. However, if both
$A$ and $B$ are finite, then $A\otimes B$ is a finite \jzs, thus a
lattice, and the following \emph{Isomorphism Theorem} holds, see
\cite{GLQu81}: 
 \begin{equation}\label{Eq:ConIsFla}
 \Conc(A\otimes B)\cong\Conc A\otimes\Conc B.
 \end{equation}
The question whether the formula \eqref{Eq:ConIsFla} extends to
the dimension monoid, \emph{i.e.}, whether the following formula
holds 
 \begin{equation}\label{Eq:SimIsFla}
 \Dim(A\otimes B)\cong\Dim A\otimes\Dim B
 \end{equation}
(the $\otimes$ on the right hand side of \eqref{Eq:SimIsFla} is
the \emph{tensor product of \cm s}, see P.A. Grillet \cite{Gril} or
F. Wehrung \cite{Wehr96})
is thus quite natural. Unfortunately, this is not the case in
general, for example, for $A=B=M_3$, the five element modular
nondistributive lattice, $A\otimes B$ is simple and not modular,
whence $\Dim(A\otimes B)$ is isomorphic to $\two$, the two-element
semilattice. However, $\Dim A=\Dim B\cong\ZZ^+$, whence
$\Dim A\otimes\Dim B\cong\ZZ^+$ again. The reason for this problem is
that modularity is not preserved under tensor product. We shall now
see how this problem can be solved for \jsd\ lattices, thus making it
possible to prove a variant of \eqref{Eq:SimIsFla} for those lattices.
We first prove a very general result.

\begin{proposition}\label{P:CanDimHom}
Let $A$ and $B$ be lattices
with zero, let $C$ be a subset of $A\ootimes B$ that satisfies the
following properties:
\begin{enumerate}
\item $(C,\subseteq)$ is a lattice;

\item $C$ is closed under finite intersection;

\item $C$ contains as elements all the mixed tensors.
\end{enumerate}

Then there exists a unique monoid homomorphism $\pi\colon\Dim
A\otimes\Dim B\to\Dim C$ such the formula 
 \[
 \pi(\DD_A(a,a')\otimes\DD_B(b,b'))=
 \DD_C((a\otimes b')\cup(a'\otimes b),a'\otimes b')
 \]
holds for all $a\leq a'$ in $A$ and $b\leq b'$ in $B$.
\end{proposition}

We observe that if $A\otimes B$ is a lattice, then it obviously
satisfies the conditions (i)--(iii) above.

\begin{proof}
We first fix $a\leq a'$ in $A$. Let
$f_{a,a'}\colon\diag{B}\to\Dim C$ be the map defined by the rule
 \[
 f_{a,a'}(x,y)=\DD_C((a\otimes y)\cup(a'\otimes x),a'\otimes y),
 \text{ for all }(x,y)\in\diag{B}.
 \]
{}From $a\leq a'$ follows that $f_{a,a'}(x,x)=0$ for all $x\in
B$. 

Now let $x\leq y\leq z$ in $B$. {}From the easily verified
relation (that holds in $C$)
 \[
 [(a\otimes y)\cup(a'\otimes x),a'\otimes y]\nearrow
 [(a\otimes z)\cup(a'\otimes x),(a\otimes z)\cup(a'\otimes y)]
 \]
follows that
 \begin{align*}
 f_{a,a'}(x,y)&+f_{a,a'}(y,z)=
 \DD_C((a\otimes y)\cup(a'\otimes x),a'\otimes y)+ \DD_C((a\otimes
 z)\cup(a'\otimes y),a'\otimes z)\\ &=\DD_C((a\otimes
 z)\cup(a'\otimes x),(a\otimes z)\cup(a'\otimes y))
 +\DD_C((a\otimes z)\cup(a'\otimes y),a'\otimes z)\\
 &=\DD_C((a\otimes z)\cup(a'\otimes x),a'\otimes z)\\
 &=f_{a,a'}(x,z).
 \end{align*}
Let $x$, $y\in B$. {}From the easily verified
relation (that holds in $C$)
 \[
 [(a\otimes x)\cup(a'\otimes(x\wedge y)),a'\otimes x]\nearrow
 [(a\otimes(x\vee y))\cup(a'\otimes y),a'\otimes(x\vee y)]
 \]
follows that
 \begin{align*}
 f_{a,a'}(x\wedge y,x)&=
 \DD_C((a\otimes x)\cup(a'\otimes(x\wedge y)),a'\otimes x)\\
 &=\DD_C((a\otimes(x\vee y))\cup(a'\otimes y),a'\otimes(x\vee
 y))\\ &=f_{a,a'}(y,x\vee y).
 \end{align*}
Therefore, we have verified that $f_{a,a'}$ is a
$\Dim C$-valued dimension function on $B$. Hence there exists a
monoid homomorphism $\varphi_{a,a'}\colon\Dim B\to\Dim C$ such
that 
 \begin{equation}\label{Eq:phiaa'xy}
 \varphi_{a,a'}(\DD_B(x,y))=f_{a,a'}(x,y),\text{ for all
 }(x,y)\in\diag{B}.
 \end{equation}
Symmetrically, for all $b\leq b'$
in $B$, we define a map $g_{b,b'}\colon\diag{A}\to\Dim C$ by the
rule 
 \[
 g_{b,b'}(x,y)=\DD_C((x\otimes b')\cup(y\otimes b),y\otimes b'),
 \text{ for all }(x,y)\in\diag{A},
 \]
and there exists a monoid homomorphism
$\psi_{b,b'}\colon\Dim A\to\Dim C$ such that 
 \begin{equation}\label{Eq:psibb'xy}
 \psi_{b,b'}(\DD_A(x,y))=g_{b,b'}(x,y),
 \text{ for all }(x,y)\in\diag{A}.
 \end{equation}
Furthermore, the symbols
$\varphi$ and $\psi$ are related as follows: 
 \begin{equation}\label{Eq:phipsi}
 \varphi_{a,a'}(\DD_B(b,b'))=\psi_{b,b'}(\DD_A(a,a'))
 =\DD_C((a\otimes b')\cup(a'\otimes b),a'\otimes b'),
 \end{equation}
for all $a\leq a'$ in $A$ and all $b\leq b'$ in
$B$. 

Now fix $b\leq b'$ in $B$. {}From \eqref{Eq:phipsi} follows that
$\varphi_{x,x}(\DD_B(b,b'))=0$ for all $x\in A$. Let $x\leq y\leq
z$ in $A$. It follows from \eqref{Eq:phipsi} and the fact that all
maps of the form either $\varphi_{u,v}$ or $\psi_{u,v}$ are monoid
homomorphisms that 
 \begin{align*}
 \varphi_{x,y}(\DD_B(b,b'))+\varphi_{y,z}(\DD_B(b,b'))&=
 \psi_{b,b'}(\DD_A(x,y)+\DD_A(y,z))\\ &=\psi_{b,b'}(\DD_A(x,z))\\
 &=\varphi_{x,z}(\DD_B(b,b')).
 \end{align*}
Similarly, for all $x$, $y\in A$,
 \begin{align*}
 \varphi_{x\wedge y,x}(\DD_B(b,b'))&=\psi_{b,b'}(\DD_A(x\wedge y,x))\\
 &=\psi_{b,b'}(\DD_A(y,x\vee y))\\
 &=\varphi_{y,x\vee y}(\DD_B(b,b')).
 \end{align*}
Therefore, for any $\beta\in\Dim B$, the map $\diag{A}\to\Dim C$,
$(x,y)\mapsto\varphi_{x,y}(\beta)$ is a dimension function on $A$, so
there exists a map $\tau_\beta\colon\Dim A\to\Dim C$ such that 
 \begin{equation}\label{Eq:taubeta}
 \tau_\beta(\DD_A(a,a'))=\varphi_{a,a'}(\beta),
 \text{ for all }a\leq a'\text{ in }A.
 \end{equation}
We define $\tau\colon\Dim A\times\Dim B\to\Dim C$
by the rule 
 \begin{equation}\label{Eq:Deftau}
 \tau(\alpha,\beta)=\tau_\beta(\alpha),
 \text{ for all }(\alpha,\beta)\in\Dim A\times\Dim B.
 \end{equation}
It follows from \eqref{Eq:Deftau} that $\tau$ is biadditive in
$\alpha$. It is also biadditive in $\beta$, because, for all
$a\leq a'$ in $A$ and all $\beta$, $\gamma\in\Dim B$,
 \begin{align*}
 \tau_{\beta+\gamma}(\DD_A(a,a'))&=\varphi_{a,a'}(\beta+\gamma)\\
 &=\varphi_{a,a'}(\beta)+\varphi_{a,a'}(\gamma)\\
 &=\tau_\beta(\DD_A(a,a'))+\tau_\gamma(\DD_A(a,a')),
 \end{align*}
whence $\tau_{\beta+\gamma}=\tau_\beta+\tau_\gamma$. Moreover,
$\tau_\beta(0)=0$ for all $\beta\in\Dim B$ and it follows from
\eqref{Eq:taubeta} that $\tau_0(\alpha)=0$ for all
$\alpha\in\Dim A$. 

Hence there exists a unique monoid homomorphism $\pi\colon\Dim
A\otimes\Dim B\to\Dim C$ such that
$\pi(\alpha\otimes\beta)=\tau(\alpha,\beta)$ for all
$(\alpha,\beta)\in\Dim A\times\Dim B$. For $a\leq a'$ in $A$ and
$b\leq b'$ in $B$,
 \begin{align*}
 \pi(\DD_A(a,a')\otimes\DD_B(b,b'))&=\tau_{\DD_B(b,b')}(\DD_A(a,a'))\\
 &=\varphi_{a,a'}(\DD_B(b,b'))\\
 &=\DD_C((a\otimes b')\cup(a'\otimes b),a'\otimes b'),
 \end{align*}
thus $\pi$ is as
required. Since the elements of the form
$\DD_A(a,a')\otimes\DD_B(b,b')$ are generators of the monoid $\Dim
A\otimes\Dim B$, the uniqueness statement is obvious.
\end{proof}

For $A$ and $B$ finite, the lattices $C$ that satisfy the
conditions of Proposition~\textup{\ref{P:CanDimHom}} are called
\emph{sub-tensor products of $A$ and $B$} in \cite{GrWe2}.

\begin{lemma}\label{L:piSurj}
Let $A$ and $B$ be finite lattices,
let $C$ be a sub-tensor product of $A$ and~$B$. Then the following
assertions hold:
\begin{enumerate}
\item $\J(C)=\setm{a\otimes b}{(a,b)\in\J(A)\times\J(B)}$.
Furthermore, the equality $(a\otimes b)_*=(a_*\otimes
b)\cup(a\otimes b_*)$ holds for all $(a,b)\in\J(A)\times\J(B)$.

\item The canonical map $\pi\colon\Dim A\otimes\Dim B\to\Dim C$ is
surjective.
\end{enumerate}
\end{lemma}

\begin{proof}
(i) Every element of $C$ is a finite join of
elements of the form $a\otimes b$ for $(a,b)\in\J(A)\times\J(B)$,
thus we obtain 
 \[
 \J(C)\subseteq\setm{a\otimes b}{(a,b)\in\J(A)\times\J(B)}.
 \]
Conversely, we put $U=(a_*\otimes b)\cup(a\otimes b_*)$, so $U\in
C$ and $U<a\otimes b$. Let $H\subset a\otimes b$ be an element of
$C$. Suppose that $H\not\subseteq U$. Then there exists $(x,y)\in
H\setminus U$. Hence $0_A<x\leq a$, $0_B<y\leq b$, and $x\nleq
a_*$ and $y\nleq b_*$, whence $x=a$ and $y=b$, so $(a,b)\in H$, a
contradiction. 

(ii) For all $X$, $Y\in C$ such that $X\prec Y$, there exists,
by Lemma~\ref{L:manylessdot}, $P\in\J(C)$ such that
$[P_*,P]\nearrow[X,Y]$. Furthermore, it follows from (i) that
there exists $(a,b)\in\J(A)\times\J(B)$ such that $P=a\otimes b$,
whence
$\DD_C(X,Y)=\DD_C(P_*,P)=\pi(\DD_A(a_*,a)\otimes\DD_B(b_*,b))$
belongs to the range of $\pi$. For the general case where $X\leq
Y$, there are $n<\omega$ and a chain $X=Z_0\prec
Z_1\prec\cdots\prec Z_n=Y$, whence
$\DD_C(X,Y)=\sum_{i<n}\DD_C(Z_i,Z_{i+1})$ belongs to the range of
$\pi$ again. Therefore, $\pi$ is surjective.
\end{proof}

Now a simple lemma about tensor products of \jz-semilattices: 

\begin{lemma}\label{L:TPSem}
Let $S$ and $T$ be \jz-semilattices,
let $a$, $a'\in S\setminus\{0_S\}$, let $b$, $b'\in
T\setminus\{0_T\}$. Then $a\otimes b\leq a'\otimes b'$ if{f}
$a\leq a'$ and $b\leq b'$.
\end{lemma}

\begin{proof}
This follows immediately from the representation of
$S\otimes T$ as the lattice of bi-ideals of $S\times T$, see
\cite{GrWe2}.
\end{proof}

Next, we recall a statement from \cite{GrWe2}: 

\begin{lemma}\label{L:Embstp}
Let $A$ and $B$ be lattices with
zero, let $C$ be a sub-tensor product of $A$ and $B$. Then there
exists a unique embedding of \jz-semilattices\newline
$\eps\colon\Conc A\otimes\Conc B\into\Conc C$ such that 
 \begin{equation}\label{Eq:aa'otbb'}
 \eps\bigl(\Theta_A(a,a')\otimes\Theta_B(b,b')\bigr)=
 \Theta_C((a\otimes b')\vee(a'\otimes b),a'\otimes b')
 \end{equation}
holds, for all $a\leq a'$ in $A$ and $b\leq b'$ in
$B$.
\end{lemma}

\begin{proof}
By \cite[Proposition~5.1 and Lemma~5.3]{GrWe2},
there exists a (necessarily unique) \jz-homomorphism
$\eps\colon\Conc A\otimes\Conc B\to\Conc C$ such that
\eqref{Eq:aa'otbb'} holds for all $a\leq a'$ in $A$ and $b\leq b'$
in $B$. By \cite[Theorem~1]{GrWe2}, $\eps$ is an embedding.
\end{proof}

By combining these results and others of this paper, we thus
obtain the following:

\begin{theorem}\label{T:DimTens}
Let $A$ and $B$ be finite \jsd\ lattices, let $C$ be a \jsd\ sub-tensor
product of $A$ and $B$. Then the canonical map $\pi\colon\Dim
A\otimes\Dim B\to\Dim C$ is an isomorphism.
\end{theorem}

\begin{proof}
It follows from Lemma~\ref{L:piSurj}(i) and
Corollary~\ref{C:DimBCF} that $\Dim C$ is the \cm\ defined by
generators $\DD(a\otimes b)$, for $(a,b)\in\J(A)\times\J(B)$, and
relations 
 \[
 \DD(a\otimes b)\ll\DD(a'\otimes b'),\text{ for all }a,\,a'\in\J(A)
 \text{ and }b,b'\in\J(B)\text{ such that }
 (a\otimes b)\DR(a'\otimes b').
 \]
Furthermore, the map $\pi$ sends $\DD_A(a)\otimes\DD_B(b)$ to
$\DD_C(a\otimes b)$, for all $(a,b)\in\J(A)\times\J(B)$. Hence, it
suffices to prove that $(a\otimes b)\DR(a'\otimes b')$ implies
that $\DD_A(a)\otimes\DD_B(b)\ll\DD_A(a')\otimes\DD_B(b')$, for
all $a$, $a'\in\J(A)$ and $b$, $b'\in\J(B)$.

By the easy direction of \cite[Lemma~2.36]{FJN}, the condition
$(a\otimes b)\DR(a'\otimes b')$ implies that $\Theta_C(a\otimes
b)\subseteq\Theta_C(a'\otimes b')$, where, for every $P\in\J(C)$,
we put $\Theta_C(P)=\Theta_C(P_*,P)$, the principal congruence of
$C$ generated by the pair $(P_*,P)$. But for $P=u\otimes v$ where
$(u,v)\in\J(A)\times\J(B)$, it follows from Lemma~\ref{L:Embstp}
that $\Theta_C(u\otimes v)=\eps(\Theta_A(u)\otimes\Theta_B(v))$,
where $\eps$ is the canonical isomorphism from $\Conc
A\otimes\Conc B$ onto $\Conc C$. Therefore, we have obtained that
 \[
 \Theta_A(a)\otimes\Theta_B(b)\leq
 \Theta_A(a')\otimes\Theta_B(b').
 \]
Therefore, by Lemma~\ref{L:TPSem},
$\Theta_A(a)\subseteq\Theta_A(a')$ and
$\Theta_B(b)\subseteq\Theta_B(b')$, whence, by
\cite[Lemma~2.36]{FJN}, $a\utr_Aa'$ and $b\utr_Bb'$. We cannot
have simultaneously $a=a'$ and $b=b'$, otherwise $a\otimes b=a'\otimes
b'$; thus either $a\tr_Aa'$ or $b\tr_Bb'$. Therefore, by
Corollary~\ref{C:tr2ll}(ii), either $\DD_A(a)=\DD_A(a')$ or
$\DD_A(a)\ll\DD_A(a')$, and either $\DD_B(b)=\DD_B(b')$ or
$\DD_B(b)\ll\DD_B(b')$, with at least one occurrence of $\ll$
taking place. This obviously implies that
$\DD_A(a)\otimes\DD_B(b)\ll\DD_A(a')\otimes\DD_B(b')$, which is
the desired conclusion.
\end{proof}

There is still the nontrivial problem left whether the statement
of Theorem~\ref{T:DimTens} is not vacuous, \emph{i.e.}, whether
for finite \jsd\ lattices $A$ and $B$, there exists a \jsd\
sub-tensor product of $A$ and $B$. We shall answer this question
affirmatively, and discuss it further, in the coming sections.

\section{Box products of \jsd\ lattices}\label{S:BPSD+} 

We start with an easy lemma, that slightly generalizes 
\cite[Lemma~1.2]{AGT}:

\begin{lemma}\label{L:G-G+}
Let $L$ be a lattice, let $G_+$ and
$G_-$ be subsets of $L$ such that every element of $L$ is a finite
join (resp., a finite meet) of elements of $G_+$ (resp., of
$G_-$). We assume that
 \[
 a\vee b=a\vee c\Rightarrow a\vee b=a\vee(b\wedge c),\quad
 \text{for all }a\in G_-\text{ and }b,\,c\in G_+.
 \]
Then $L$ is \jsd.
\end{lemma}

\begin{proof}
By \cite[Lemma~1.2]{AGT}, it suffices to prove
that $a\vee b=a\vee c$ implies that $a\vee b=a\vee(b\wedge c)$,
for all $a\in L$ and all $b$, $c\in G_+$. Put $d=a\vee b=a\vee c$,
and suppose that $a\vee(b\wedge c)<d$. By assumption on
$G_-$, there are $n>0$ and $e_0$, \dots, $e_{n-1}\in G_-$ such
that $a\vee(b\wedge c)=\bigwedge_{i<n}e_i$. Hence there exists
$i<n$ such that $d\nleq e_i$. Moreover,
 \[
 e_i\vee b=e_i\vee a\vee(b\wedge c)\vee b=e_i\vee d=e_i\vee c,
 \]
with $e_i\in G_-$ and $b$, $c\in G_+$, therefore, by assumption, 
 \[
 e_i\vee b=e_i\vee(b\wedge c)=e_i,
 \]
whence $b\leq e_i$. Thus, $d=a\vee b\leq e_i$, a contradiction.
\end{proof}

For lattices $A$ and $B$ both with least and greatest element, the
\emph{box product} $A\bp B$ of $A$ and $B$ is a particular case of
sub-tensor product of $A$ and $B$, see \cite{GrWe4}. In this case,
the elements of $A\bp B$ are exactly the finite intersections of
\emph{mixed tensors} defined in Section~\ref{S:DimTens}.

\begin{corollary}\label{C:BPSD+}
For any \jsd\ lattices $A$ and
$B$, the box product $A\bp B$ is \jsd.
\end{corollary}

\begin{proof}
We use the notation and terminology of \cite{GrWe4}.
Write $A=\varinjlim_{a\in A}\upw a$ and $B=\varinjlim_{b\in B}\upw b$,
with the obvious transition homomorphisms and limiting maps. Then
$A\bp B=\varinjlim_{(a,b)\in A\times B}(\upw a\bp\upw b)$, with all
the lattices $\upw a$ (for $a\in A$) and $\upw b$ (for $b\in B$) \jsd,
thus it suffices to consider the case where both $A$ and $B$ are
lattices with zero. Next, if $A'$ (resp., $B'$) is the lattice
obtained by adding a new unit to $A$ (resp., $B$), then both $A'$ and
$B'$ are \jsd\ and $A\bp B$ is isomorphic to an ideal of $A'\bp B'$.
Therefore, we have reduced the problem to \emph{bounded} lattices $A$
and $B$. 

By the definition of the box product, the set $G_-$ defined by 
 \[
 G_-=\setm{a\bp b}{(a,b)\in A\times B}
 \]
generates $(A\bp B,\wedge)$, while, since both $A$ and $B$ are
bounded, the set $G_+$ defined by
 \[
 G_+=\setm{a\ltp b}{(a,b)\in A\times B}
 \]
generates $(A\bp B,\vee)$. Hence, by Lemma~\ref{L:G-G+}, to
verify that $A\bp B$ is \jsd, it suffices to verify that 
 \begin{equation}\label{Eq:HypSD}
 (a\bp b)\vee(x_0\ltp y_0)=(a\bp b)\vee(x_1\ltp y_1)
 \end{equation}
implies that
 \begin{equation}\label{Eq:ConcSD}
 (a\bp b)\vee(x_0\ltp y_0)=(a\bp b)\vee(x\ltp y),
 \end{equation}
where we put
 \[
 x=x_0\wedge x_1\text{ and }y=y_0\wedge y_1
 \]
(we use the fact
that $(x_0\ltp y_0)\cap(x_1\ltp y_1)=x\ltp y$). The conclusion is
trivial if $x_0\ltp y_0\leq a\bp b$, so suppose that $x_0\ltp
y_0\nleq a\bp b$, and thus also $x_1\ltp y_1\nleq a\bp b$, so that 
 \begin{equation}\label{Eq:Nleq}
 x_i\nleq a\text{ and }y_i\nleq b,\text{ for all }i<2.
 \end{equation}
Furthermore, it is easy to verify that
 \begin{equation}\label{Eq:ExpBP}
 \begin{aligned}
 (a\bp b)\vee(u\ltp v)&=((a\vee u)\bp b)\cap(a\bp(b\vee v))\\
 &=(a\bp b)\cup((a\vee u)\circ(b\vee v)),
 \text{ for all }(u,v)\in A\times B.
 \end{aligned}
 \end{equation}
Therefore, by using \eqref{Eq:HypSD},
\eqref{Eq:Nleq}, and \eqref{Eq:ExpBP}, we obtain that
 \[
 a\vee x_0=a\vee x_1\text{ and }b\vee y_0=b\vee y_1,
 \]
from which it follows, since both $A$ and $B$ are \jsd, that 
\[
 a\vee x_0=a\vee x\text{ and }b\vee y_0=b\vee y,
 \]
so applying \eqref{Eq:ExpBP} to the pairs $(x_0,y_0)$ and $(x,y)$
yields the conclusion \eqref{Eq:ConcSD}.
\end{proof}

As an immediate consequence of
Theorem~\ref{T:DimTens} and Corollary~\ref{C:BPSD+}, we observe the
following:

\begin{corollary}\label{C:AbpB}
Let $A$ and $B$ be finite \jsd\ lattices. Then the relation
$\Dim(A\bp B)\cong\Dim A\otimes\Dim B$ holds.
\end{corollary}

Another result related to Corollary~\ref{C:BPSD+} is the following: 

\begin{proposition}\label{P:LBSD+}
Let $A$ and $L$ be lattices
with zero, with $A$ finite. If $A$ is a lower bounded homomorphic
image of a free lattice and if $L$ is \jsd, then $A\otimes L$ is \jsd.
\end{proposition}

\begin{note}
An immediate application of
\cite[Proposition 2.9]{GrWe2} and \cite[Corollary~5.4]{GrWe3} yields
that \emph{if $A$ and $B$ are finite lower bounded homomorphic images
of free lattices, then $A\otimes B$ is a lower bounded homomorphic
image of a free lattice}.
\end{note}

\begin{proof}
Put $P=\J(A)$. For any $p\in P$, we put
 \[
 \CM(p)=\setm{I\subseteq P}
 {I\text{ is a minimal nontrivial join-cover of }p}.
 \]
We recall that the join-dependency relation
$\DR$ on $A$ can be defined by 
 \[
 p\DR q\text{ if and only if }\exists I\in\CM(p)
 \text{ such that }q\in I.
 \]
Let $x\colon P\to L$ be an antitone map. The
\emph{adjustment sequence} of $x$ is defined by $x^{(0)}=x$, and
$x^{(n+1)}=(x^{(n)})^{(1)}$, where $x^{(1)}$ is defined by the rule
 \[
 x^{(1)}(p)=x(p)\vee\bigvee_{I\in\CM(p)}\bigwedge_{q\in I}x(q),
 \text{ for all }p\in P.
 \]
(In particular, the map $x^{(n)}$, for any $n\in\omega$, is
still antitone.) By \cite[Remark~6.6 and Theorem~4(iii)]{GrWe3},
since $A$ is a finite lower bounded homomorphic image of a free lattice
(`amenable'), the adjustment sequence of any antitone map
$x\colon P\to L$ is eventually constant, hence $A\otimes L\cong
A[L]$, where $A[L]$ is defined as the set of all antitone maps
$x\colon P\to L$ such that 
 \[
 \bigwedge_{q\in I}x(q)\leq x(p),\text{ for every }p\in P
 \text{ and every minimal nontrivial join-cover }I\text{ of }p.
 \]
Let $y'$ (resp., $z'$) be the antitone maps from $P$ to $L$ defined
by the rules
\[
 y'(p)=x(p)\vee y(p)\text{ and }z'(p)=x(p)\vee z(p),
 \text{ for all }p\in P.
 \]
 Hence $x\vee y$ (resp., $x\vee z$) is the supremum
of the [eventually constant] adjustment sequence of $y'$ (resp.,
$z'$). 

To conclude the proof, it suffices to prove that $A[L]$ is \jsd.
So, let $x$, $y$, $z\in A[L]$ such that $x\vee y=x\vee z$ and
$y\wedge z\leq x$, we prove that $y\leq x$ (and so also $z\leq
x$). For this, we prove, by downward $\DR$-induction, that the
following equality holds for all $p\in P$:
 \begin{equation}\label{Eq:xvyvzx}
 (x\vee y)(p)=(x\vee z)(p)=x(p).
 \end{equation}
Suppose that \eqref{Eq:xvyvzx} holds for all $q\in
P$ such that $p\DR q$. In particular,
$(y')^{(n)}(q)=(z')^{(n)}(q)=x(q)$ for any $n\in\omega$ and any
$q\in P$ such that $p\DR q$. We prove that \eqref{Eq:xvyvzx} holds
at $p$. To achieve this, we first prove a claim.

\goodbreak
\setcounter{claim}{0}
\begin{claim}\label{Cl:xyz}
$(y')^{(n)}(p)=y'(p)$ and $(z')^{(n)}(p)=z'(p)$, for all
$n\in\omega$.
\end{claim}

\begin{cproof}
We argue by induction on $n$. The result is trivial
for $n=0$. Suppose that it holds for $n$. For any $I\in\CM(p)$ and
any $q\in I$, the relation $p\DR q$ holds, thus, by the induction
hypothesis, $(x\vee y)(q)=x(q)$, whence $(y')^{(n)}(q)=x(q)$.
Therefore, we can compute
 \begin{align*}
 (y')^{(n+1)}(p)&=
 (y')^{(n)}(p)\vee\bigvee_{I\in\CM(p)}
 \bigwedge_{q\in I}(y')^{(n)}(q)\\
 &=y'(p)\vee\bigvee_{I\in\CM(p)}\bigwedge_{q\in I}x(q)&&
 \quad\text{(by the induction hypotheses)}\\
 &=\biggl(x(p)\vee\bigvee_{I\in\CM(p)}\bigwedge_{q\in I}
 x(q)\biggr)\vee y(p)\\ &=x(p)\vee y(p)&&
 \quad(\text{since }x\in A[L])\\
 &=y'(p).
 \end{align*}
Similarly, we can prove that $(z')^{(n)}(p)=z'(p)$,
for all $n\in\omega$.
\end{cproof}

As an immediate consequence of Claim~\ref{Cl:xyz}, the two
following equalities hold:
 \begin{equation}\label{Eq:xvyz}
 (x\vee y)(p)=x(p)\vee y(p)\text{ and }(x\vee z)(p)=x(p)\vee z(p).
 \end{equation}
Hence the
assumption that $x\vee y=x\vee z$, plus the fact that $y(p)\wedge
z(p)=(y\wedge z)(p)\leq x(p)$ and the \jsdy\ of $L$, imply that
$y(p),z(p)\leq x(p)$, which, again by \eqref{Eq:xvyz}, implies
that $(x\vee y)(p)=(x\vee z)(p)=x(p)$.

So we have established that \eqref{Eq:xvyvzx} holds for every
$p\in P$, whence $x\vee y=x\vee z=x$. Therefore, $y,z\leq x$,
which completes the proof.
\end{proof}

The following section shows that tensor products are not as
well-behaved, for finite \jsd\ lattices, as box products. 

\section{Non-preservation of join-semidistributivity by tensor
product}\label{S:TsCtEx}

Let $L$ be the lattice of all order-convex subsets of a
four-element chain, see Figure~2. Hence $L$ is atomistic and \jsd,
and $\J(L)=\{a,b,a',b'\}$, with the generators $a$,
$b$, $a'$, $b'$ subjected to the following relations: 
 \begin{align}
 a&<a'\vee b;\label{Eq:aa'b}\\
 a'&<a\vee b';\label{Eq:a'ab}\\
 a,a'&<b\vee b'.\label{Eq:aa'bb'}
 \end{align}

\begin{figure}[hbt]
\includegraphics{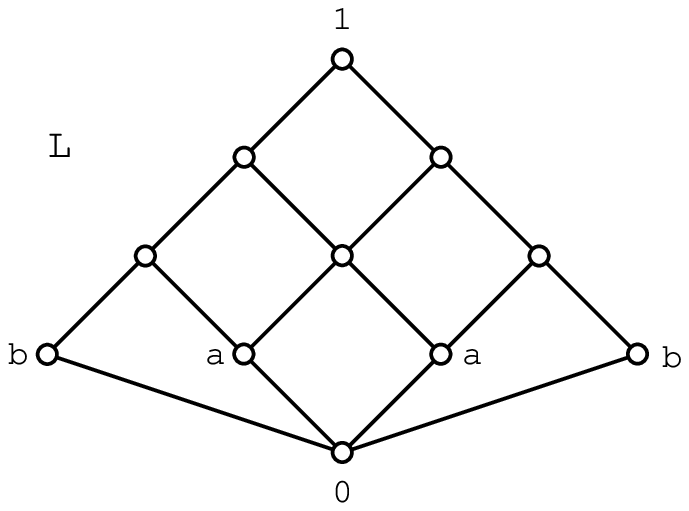}
\caption{}
\end{figure}

We shall now prove that $L\otimes L$ is not \jsd. We define an
element $H$ of $L\otimes L$ as follows:
 \[
 H=(a\otimes b)\vee(b\otimes a')\vee
 (a'\otimes b')\vee(b'\otimes a).
 \]
We use the representation of the tensor product $L\otimes L$ as the
lattice of bi-ideals of $L\times L$. Since $a$, $b$,
$a'$, $b'$ are distinct atoms of $L$, we also have
 \begin{equation}\label{Eq:capped}
 H=(a\otimes b)\cup(b\otimes a')\cup(a'\otimes b')\cup(b'\otimes a).
 \end{equation}
(Indeed, it
suffices to verify that the right hand side of \eqref{Eq:capped}
is a bi-ideal of $L\times L$.) In particular, we obtain that 
 \begin{equation}\label{Eq:nin}
 a\otimes a\nleq H.
 \end{equation}
Furthermore, by using \eqref{Eq:aa'b}, we obtain
the inequalities
 \begin{align*}
 a\otimes a&\leq(a\otimes b)\vee(a\otimes a'),\\
 a\otimes a'&\leq(b\otimes a')\vee(a'\otimes a'),
 \end{align*}
from
which it follows that
 \begin{align}
 a\otimes a&\leq(a\otimes a')\vee H,\label{Eq:aaa'H}\\
 a\otimes a'&\leq(a'\otimes a')\vee H.\label{Eq:aa'a'H}
 \end{align}
Similarly, by using \eqref{Eq:a'ab}, we obtain the inequalities 
 \begin{align*}
 a'\otimes a'&\leq(a'\otimes a)\vee(a'\otimes b'),\\
 a'\otimes a&\leq(a\otimes a)\vee(b'\otimes a),
 \end{align*}
from which it follows that
 \begin{align}
 a'\otimes a'&\leq(a'\otimes a)\vee H,\label{Eq:a'a'b'H}\\
 a'\otimes a&\leq(a\otimes a)\vee H.\label{Eq:a'ab'H}
 \end{align}
{F}rom the inequalities
\eqref{Eq:aaa'H}, \eqref{Eq:aa'a'H}, \eqref{Eq:a'a'b'H}, and
\eqref{Eq:a'ab'H} follows that 
 \[
 (a\otimes a)\vee H=(a\otimes a')\vee H=(a'\otimes a')\vee
 H=(a'\otimes a)\vee H.
 \]
By \eqref{Eq:nin} and since
$(a\otimes a)\wedge(a\otimes a')=0$, it follows that $L\otimes L$
is not \jsd.

\section{Open problems}\label{S:OpenPbs}

Our first problem is motivated by the so-called \emph{Separativity
Conjecture} in ring theory, that asks, for, say, a von Neumann
regular ring $R$, whether the monoid $V(R)$ of all isomorphism
classes of finitely generated projective right $R$-modules is
separative:

\begin{problem}\label{Pb:SepConj}
For a lattice $L$, is $\Dim L$ separative?
\end{problem}

By the results of \cite{WDim}, if the Separativity Conjecture fails
for rings, then it also fails for lattices, and even for complemented
modular lattices. The converse is not clear, although it may
shed some light on the ring theoretical problem. An equivalent form
of Problem~\ref{Pb:SepConj} is to ask, for a natural number $n$,
whether $\Dim\FL(n)$ is separative, where $\FL(n)$ denotes the free
lattice on $n$ generators.

\begin{problem}\label{Pb:DimPres}
Let $K$ be a finite \jsd\
lattice. Can $K$ be embedded into a finite, atomistic, \jsd\
lattice $L$ in a \emph{dimension-preserving way}, \emph{i.e.}, in
such a way that the canonical map from $\Dim K$ to $\Dim L$ is an
isomorphism?
\end{problem}


For $K$ a lower bounded homomorphic image of a free lattice, one can
prove, by using methods from this paper, that a positive solution to
Problem~\ref{Pb:DimPres} is provided by Tischendorf's extension, see
M. Tischendorf \cite{Tisch}. For different classes of lattices $K$ we
cannot hope a positive solution of Problem~\ref{Pb:DimPres},
\emph{e.g.}, let $K$ be a finite modular lattice that cannot be
embedded into any finite atomistic modular lattice (the subgroup
lattice of $(\ZZ/4\ZZ)^3$ is such an example, see C.~Herrmann and
A.P.~Huhn \cite{HeHu}). Then $K$ cannot be embedded
dimension-preservingly into any finite atomistic lattice $L$, for
$\Dim L$ is isomorphic to $\Dim K$, thus it is cancellative, thus
$L$ is modular.

On the other hand, it has been proved that every finite \jsd\
lattice can be embedded into a finite atomistic \jsd\ lattice, see
\cite[Theorem~1.11]{AGT}.

Our next problem calls for a generalization of
Theorem~\ref{T:DimTens} to the infinite case:

\begin{problem}\label{Pb:AbpB}
Let $A$ and $B$ be bounded lattices. Is the canonical map\linebreak
$\pi\colon\Dim A\otimes\Dim B\to\Dim(A\bp B)$ (see
Proposition~\ref{P:CanDimHom}) surjective? If $A$ is \jsd, is $\pi$
an isomorphism?
\end{problem}

If both $A$ and $B$ are finite \jsd, then the required isomorphy
holds, see Corollary~\ref{C:AbpB}.

\begin{problem}\label{Pb:CharCon}
Characterize the dimension monoids of finite lattices (\jsd\ finite
lattices, finite lower bounded homomorphic images of free lattices,
respectively).
\end{problem}

Our next problem is motivated by the result of
Section~\ref{S:TsCtEx}:

\begin{problem}\label{Pb:LotL}
Let $L$ be a finite lattice. If $L\otimes L$ is \jsd, is $L$ a lower
bounded homomorphic image of a free lattice?
\end{problem}

\begin{problem}
Does there exist for the dimension theory of \jsd\ lattices an
analogue of the continuous geometries with nondiscrete dimension
monoid?
\end{problem}

For the \jsd\ lattices of the present paper, the dimension functions
are $\ZZ^+\cup\set{\infty}$-valued, while for the hypothetical new
objects, the dimension functions would be
$\RR^+\cup\set{\infty}$-valued. Can one cultivate any analogy with the
theory underlying the decomposition of, say, a self-injective von
Neumann regular ring into factors of type I, II, or III (see, for
example, K.R. Goodearl \cite{Good91})? Then one could say that the
results of the present paper deal essentially with type I, although
the relation $\DR^\infty$ introduced in Section~\ref{S:JDD}
definitely carries a touch of type III.

\end{document}